\documentclass[11pt]{amsart}
\usepackage{amssymb,amsmath,txfonts,mathrsfs,mathtools,titletoc, tikz, stmaryrd}
\usepackage[alphabetic]{amsrefs}

\newtheorem{theorem}{Theorem}[section]
\newtheorem{prop}[theorem]{Proposition}
\newtheorem{lemma}[theorem]{Lemma}
\newtheorem{remark}[theorem]{Remark}

\newtheorem{question}[theorem]{Question}

\newtheorem{cor}[theorem]{Corollary}

\newtheorem{conj}[theorem]{Conjecture}

\numberwithin{equation}{section}

\def\pf{{\it Proof:}~}

\begin{document}

\title[The growth rate of harmonic functions]{The growth rate of harmonic functions}
\author{Guoyi Xu}
\address{Department of Mathematical Sciences\\Tsinghua University, Beijing\\P. R. China, 100084}
\email{guoyixu@tsinghua.edu.cn}
\date{\today}
\date{\today}

\begin{abstract}
We study the growth rate of harmonic functions in two aspects: gradient estimate and frequency. We obtain the sharp gradient estimate of positive harmonic function in geodesic ball of complete surface with non-negative curvature. On complete Riemannian manifolds with non-negative Ricci curvature and maximal volume growth, further assume the dimension of the manifold is not less than three, we prove that quantitative strong unique continuation yields the existence of non-constant polynomial growth harmonic functions. Also the uniform bound of frequency for linear growth harmonic functions on such manifolds is obtained, and this confirms a special case of Colding-Minicozzi's conjecture on frequency.
\\[3mm]
\end{abstract}
\thanks{The author was partially supported by NSFC-11771230, Beijing Natural Science Foundation Z190003}

\maketitle

\titlecontents{section}[0em]{}{\hspace{.5em}}{}{\titlerule*[1pc]{.}\contentspage}
\titlecontents{subsection}[1.5em]{}{\hspace{.5em}}{}{\titlerule*[1pc]{.}\contentspage}
\tableofcontents

\section{Introduction}

This paper studies the harmonic function's growth rate on manifolds and its related application. Unless otherwise mentioned, all manifolds in this paper have non-negative Ricci curvature. By Yau's Liouville theorem, any positive harmonic functions on complete manifolds with non-negative Ricci curvature is constant function. We firstly study the positive harmonic functions defined on a unit geodesic ball. Note if we do scaling on manifolds with non-negative Ricci curvature, the unit geodesic ball case can be extended to general geodesic ball case. We get the following results characterizing the `largest' positive harmonic functions in the geodesic ball. Define $\mathcal{H}^+$ as the set of positive harmonic functions $u$ defined on $B_1(p)$ with $u(p)= 1$, where $B_1(p)\subseteq M^n$ is the geodesic ball centered at $p$ with radius $= 1$.
\begin{theorem}\label{thm positive harm inside}
{For any $s\in [0, 1)$, there exists $u_s\in \mathcal{H}^+$ such that
\begin{align}
\sup_{y\in \partial B_s(p)}u_s(y)= \max_{u\in \mathcal{H}^+} \sup_{y\in \partial B_s(p)} u(y) . \nonumber 
\end{align}
Furthermore there exists $x_s\in \partial B_1(p)$ (possibly not unique) such that $ \displaystyle u_s(y)= \frac{P(x_s, y)}{P(x_s, p)}$.  
}
\end{theorem}

Motivated by the gradient estimate of Cheng and Yau \cite{CY}, we study the largest growth rate of positive harmonic functions in a unit geodesic ball. So far we can only prove the sharp gradient estimate on surfaces.
\begin{theorem}\label{thm sharp 2-dim gradient est}
{Assume $u: B_{1}(p)\rightarrow \mathbb{R}^+$ is a harmonic function, where $M^2$ is a Riemannian manifold with $Rc\geq 0$ and $B_1(p)\subseteq M^2$ is a geodesic ball centered at $p$ with radius $= 1$, then $\sup\limits_{x\in B_s(p)}\frac{|\nabla u|}{u}\leq \frac{1}{1- s}+ \frac{1}{1+ s}$ for any $s\in [0, 1)$.
}
\end{theorem}

The model of the above theorem is the sharp gradient estimate obtained on $\mathbb{R}^n$, which is exactly the gradient estimate of the corresponding Poisson kernel on the unit ball. The sharp gradient estimate we are looking for, is one step to reveal the largest growth rate of positive harmonic function on the geodesic balls of complete manifolds with $Rc\geq 0$, which is motivated by the study of polynomial growth harmonic function's frequency. For higher dimensional manifolds (the dimension is greater than $2$), this type sharp gradient estimate is unknown so far.

It is well known that on any complete noncompact manifold, there always exist nonconstant harmonic functions (see \cite{GW}). The proof of such existence result is based on the classical unique continuation for harmonic functions. The frequency for harmonic function was firstly introduced by Almgren \cite{Almgren}. Later Garofalo and Lin \cite{GL} used the bound on frequency to prove the unique continuation, which partly revealed the deep connection between the existence problem and the bound on frequency.

The polynomial growth harmonic functions are important harmonic functions. Yau \cite{problem} conjectured: on a complete manifold with $Rc\geq 0$, the space of harmonic functions with polynomial growth is finite dimensional. In $1997$, this conjecture was completely proved by Colding and Minicozzi \cite{CMAnn} (also see \cite{Li2}, \cite{CMcpam}, \cite{CMinv} for further developments). 

Sormani \cite{Sormani} proved that the existence of any nonconstant harmonic functions with polynomial growth on linear volume growth manifolds with $Rc\geq 0$ leads to the splitting of the manifolds. This can be used to construct the manifolds with $Rc\geq 0$, which does not admit any nonconstant harmonic functions with polynomial growth. 

For a complete Riemannian manifold $M^n$ with $Rc\geq 0$, define $\mathrm{V}_{M}= \lim\limits_{r\rightarrow \infty}\frac{\mathrm{Vol}(B_r(p))}{\omega_nr^n}$; if $V_{M}> 0$, we say that $M^n$ has \textbf{maximal volume growth}. Although Yau's conjecture was confirmed, we do not know whether there always exists a non-constant polynomial growth harmonic function on $M^n$ with maximal volume growth. In the rest of the introduction, we assume $M^n$ has maximal volume growth and $n\geq 3$.

Colding and Minicozzi \cite{CM-JDG-97} introduced a slightly different frequency function $\mathscr{F}_u(r)$ for harmonic function $u$ on $M^n$, and obtained some local estimates of $\mathscr{F}_u(r)$. Furthermore,  they \cite{CM} made the following conjectures about frequency of harmonic functions.
\begin{conj}\label{conj global freq upper bound conj}  
{\begin{enumerate}
\item[(a)] If $u$ is a non-constant polynomial growth harmonic function on $M^n$, then $\mathscr{F}_u(r)$ is uniformly bounded.
\item[(b)] (Quantitative strong unique continuation). Suppose $B_{2r}(p)\subseteq M^n$ for some $p\in M^n$, and $u$ is a harmonic function on $B_{2r}(p)$ with $\sup\limits_{s\in [r, 2r]}\mathscr{F}_u(s)\leq d$, where $d> 0$ is some constant; then there exists $C= C(n, d, \mathrm{V}_M)$ such that $\sup\limits_{s\in [0, r]}\mathscr{F}_u(s)\leq C$.
\end{enumerate}
}
\end{conj}

In this paper, we prove the following theorem revealing the relation between the existence of polynomial growth harmonic function and quantitative strong unique continuation conjecture.
\begin{theorem}\label{thm quc imply existence}
{If the quantitative unique continuation holds on $M^n$, then there exists a non-constant polynomial growth harmonic function.
}
\end{theorem} 

And we also confirm part $(a)$ of Conjecture \ref{conj global freq upper bound conj} for linear growth harmonic function in the following theorem.
\begin{theorem}\label{thm uniform bound of freq}
{If $u$ is a linear growth harmonic function on $M^n$, then 
$\sup\limits_{r\geq 0}\mathscr{F}_u(r)\leq C(M^n, p)$.
}
\end{theorem}

The organization of this paper is as the following. In part I, we study the harmonic functions defined in a geodesic ball with non-negative Ricci curvature. Especially, we give the sharp gradient estimate of positive harmonic functions, by detailed study of Poisson kernel in the geodesic ball. 

In part II, we study the sufficient condition guaranteeing the existence of polynomial growth harmonic functions on the whole complete Riemannian manifold in section $4$. Then we prove any linear growth harmonic functions has bounded frequency in section $5$. The key of the proof is the asymptotic mean value equality proved by Li \cite{Li-Large}, which provides us the uniform control of $u$ near the infinity of the manifolds. Finally, we show the existence of non-constant polynomial growth harmonic functions on Perelman's example manifolds, which has maximal volume growth and non-unique tangent cones at infinity.


\section*{Part I. Positive harmonic functions in a geodesic ball}

In this part, let $G(x, y)$ be the positive Dirichlet Green's function of $B_1(p)$ and $\vec{n}_x$ is the inward unit normal vector of $\partial B_1(p)$ at $x$. Then for $n\geq 2$, we have
\begin{align}
P(x, y)= \frac{\partial G(x, y)}{\partial \vec{n}_x}, \quad \quad \quad \forall x\in \partial B_1(p), y\in B_1(p) , \nonumber 
\end{align}
where $P(x, y)$ is the Poisson kernel of $B_1(p)$.

\section{Positive harmonic function and Poisson kernel}

Let $\alpha= n(n- 2)\omega_n$, where $\omega_n$ is the volume of the unit ball $B_1(0)\subseteq \mathbb{R}^n$. 

\pf[of Theorem \ref{thm positive harm inside}]
{From Cheng-Yau's gradient estimate, we know that 
\begin{align}
\sup_{B_s(p)}\frac{|\nabla u|}{u}\leq \frac{C(n)}{1- s},  \quad \quad\quad \quad s\in [0, 1). \nonumber 
\end{align}
Hence there exists $\displaystyle h(s)\vcentcolon= \sup_{u\in \mathcal{H}^+\atop x\in B_s(p)}u(x)< \infty$. From the compactness theorem of bounded harmonic functions, there exists $u_s\in \mathcal{H}^+$ such that 
\begin{align}
\sup_{x\in B_s(p)} u_s(x)= h(s) , \quad \quad \quad \quad s\in [0, 1). \nonumber 
\end{align}

Similar to \cite[Corollary $6.15$]{ABR}, there exists a measure $d\mu_s$ on $\partial B_1(p)$ such that $\displaystyle u_s(y)= \int_{\partial B_1(p)} P(x, y)d\mu_s(x)$, where $P(x, y)$ is the Poisson kernel of $B_1(p)$. From $u_s(p)= 1$, we obtain $\displaystyle \int_{\partial B_1(p)}P(x, p)d\mu_s(x)= 1$. Then
\begin{align}
u_s(y)&= \int_{\partial B_1(p)} P(x, y)d\mu_s(x)\leq \sup_{x\in \partial B_1(p)}\frac{P(x, y)}{P(x, p)}\cdot \int_{\partial B_1(p)} P(x, p)d\mu_s(x) \nonumber \\
&= \sup_{x\in \partial B_1(p)}\frac{P(x, y)}{P(x, p)} .\nonumber  
\end{align}
Furthermore, the equality holds if and only if $d\mu_s(x)= \frac{1}{P(x, p)}\mathscr{P}_{E(y)}(x)$, where $\mathscr{P}_{E(y)}(x)$ is a probability measure on $\displaystyle E(y)= \big\{\tilde{x}\in \partial B_1(p):  \frac{P(\tilde{x}, y)}{P(\tilde{x}, p)}=\sup_{x\in \partial B_1(p)}\frac{P(x, y)}{P(x, p)}\big\}$. And we define $\psi(y)\in \partial B_1(p)$ such that $\displaystyle \frac{P(\psi(y), y)}{P(\psi(y), p)}= \sup_{x\in \partial B_1(p)}\frac{P(x, y)}{P(x, p)}$. Note that $\psi(y)$ possibly is not unique.

There exists $y_s\in \partial B_s(p)$ such that $\displaystyle u_s(y_s)= \sup_{y\in B_s(p)} u_s(y)$. Also there exists $x_s\in \partial B_1(p)$ such that $\displaystyle \frac{P(x_s, y_s)}{P(x_s, p)}= \sup_{x\in \partial B_1(p)}\frac{P(x, y_s)}{P(x, p)}$. Define $\displaystyle v_s(y)= \frac{P(x_s, y)}{P(x_s, p)}: B_1(p)\rightarrow \mathbb{R}^+$, then $v_s\in \mathcal{H}^+$ and 
\begin{align}
h(s)\geq \sup_{y\in B_s(p)}v_s(y)\geq \sup_{y\in B_s(p)}u_s(y)= h(s). \nonumber 
\end{align}
Especially, one gets $\displaystyle \sup_{x\in \partial B_1(p)}\frac{P(x, y_s)}{P(x, p)}= u_s(y_s)= h(s)$.

From the above, we get $\displaystyle d\mu_s(x)= \frac{1}{P(x, p)}\mathscr{P}_{E(y_s)}(x)$, which implies
\begin{align}
u_s(y)= \int_{\partial B_1(p)} P(x, y)d\mu_s(x)= \frac{P(\psi(y_s), y)}{P(\psi(y_s), p)}. \nonumber 
\end{align}
}
\qed

\begin{lemma}\label{lem sharp gradient estimate of Mfld}
{On $M^n$ with $n\geq 2$, assume $u$ is a positive harmonic function on $B_1(p)\subseteq M^n$, then for any $s\in [0, 1)$,
\begin{align}
\sup\limits_{y\in B_s(p)}|\nabla \ln u|(y)\leq \sup\limits_{x\in \partial B_1(p)}\sup\limits_{y\in B_s(p)}|\nabla_y \ln P(x, y)| . \nonumber 
\end{align}
}
\end{lemma}

\pf
{For any $\epsilon> 0$, let $G^\epsilon(x, y)$ be the Dirichlet Green's function of $B_{1-\epsilon}(p)$, then for any $y\in \partial B_s$, we have $u(y)= \int_{\partial B_{1-\epsilon}}u(x)P^\epsilon(x, y)dx$, where $P^\epsilon(x, y)= \frac{\partial G^\epsilon(x, y)}{\partial\vec{n}_x}$. We get
\begin{align}
\frac{|\nabla u(y)|}{u(y)}= \frac{\Big|\int_{\partial B_{1-\epsilon}}u(x)\nabla_y P^\epsilon(x, y)dx\Big|}{\int_{\partial B_{1-\epsilon}}u(x)P^\epsilon(x, y)dx}\leq \sup_{x\in \partial B_{1- \epsilon}}\frac{\Big|\nabla_y P^\epsilon(x, y)\Big|}{P^\epsilon(x, y)} \nonumber
\end{align}
Letting $\epsilon\rightarrow 0$ in the above, we have
\begin{align}
\frac{|\nabla u(y)|}{u(y)}\leq \sup_{x\in \partial B_{1}}\frac{\Big|\nabla_y P(x, y)\Big|}{P(x, y)} .\label{one crucial ineq}
\end{align}
The conclusion follows from the above.
}
\qed

Now, we have the sharp gradient estimate on $\mathbb{R}^n$.
\begin{cor}\label{lem sharp gradient estimate of Rn}
{On $\mathbb{R}^n$ with $n\geq 2$, assume $u$ is a positive harmonic function on $B_1$, then $\sup\limits_{y\in B_s}\frac{|\nabla u|}{u}(y)\leq \frac{n- 1}{1- s}+ \frac{1}{1+ s}$ for any $0\leq s< 1$. Furthermore, the equality holds if and only if $u(y)= \frac{1- |y|^2}{n\omega_n |x_0- y|^n}$ for some $x_0\in \partial B_1(0)$.
}
\end{cor}

\pf
{\textbf{Step (1)}. Recall we have the Poisson kernel $P(x, y)$ of $B_1(0)$ has the following expression:
\begin{align}
P(x, y)= \frac{\partial G(x, y)}{\partial \vec{n}_x}= \frac{1- |y|^2}{n\omega_n |x- y|^n}, \quad \quad \quad \forall x\in \partial B_1(0), y\in B_1(0) \nonumber 
\end{align}

It is straightforward to get
\begin{align}
|\nabla_y \ln P(x, y)|^2= \frac{2n- 2(n- 2)s^2}{(1- s^2)^2}+ \frac{n(n- 2)}{1+ s^2- 2xy}, \quad \quad \quad \quad \forall x\in \partial B_1, y\in \partial B_s \nonumber 
\end{align}

Hence 
\begin{align}
\sup\limits_{x\in \partial B_1\atop y\in \partial B_s}\frac{|\nabla_y P(x, y)|}{P(x, y)}&= \sup\limits_{x\in \partial B_1\atop y\in \partial B_s}\sqrt{\frac{2n- 2(n- 2)s^2}{(1- s^2)^2}+ \frac{n(n- 2)}{1+ s^2- 2xy}} \nonumber \\
&= \sqrt{\frac{2n- 2(n- 2)s^2}{(1- s^2)^2}+ \frac{n(n- 2)}{1+ s^2- 2s}}= \frac{n- 1}{1- s}+ \frac{1}{1+ s} \nonumber 
\end{align}
We conclude $\sup\limits_{x\in \partial B_1\atop y\in B_s}\frac{|\nabla_y P(x, y)|}{P(x, y)}= \frac{n- 1}{1- s}+ \frac{1}{1+ s}$ by the fact that $\frac{n- 1}{1- s}+ \frac{1}{1+ s}$ is increasing in $s$. Combining the above with Lemma \ref{lem sharp gradient estimate of Mfld}, the first conclusion follows.

\textbf{Step (2)}. If the equality holds for some $u$, then we can assume that $|\nabla \ln u(y_0)|= \frac{n- 1}{1- s}+ \frac{1}{1+ s}$ for some $y_0\in \partial B_s$. From (\ref{one crucial ineq}), we have
\begin{align}
\frac{n- 1}{1- s}+ \frac{1}{1+ s}= |\nabla \ln u(y_0)|\leq \sup_{x\in \partial B_{1}}\Big|\nabla_y \ln P(x, y_0)\Big|= \frac{n- 1}{1- s}+ \frac{1}{1+ s} \label{integral on the boundary}
\end{align}

Assume $\Big|\nabla_y \ln P(x_0, y_0)\Big|= \frac{n- 1}{1- s}+ \frac{1}{1+ s}$ for some $x_0\in \partial B_1$, then it is easy to see that $\Big|\nabla_y \ln P(x_1, y_0)\Big|< \Big|\nabla_y \ln P(x_0, y_0)\Big|$ for any $x_1\neq x_0$. Combining this inequality with (\ref{integral on the boundary}), we get $u\big|_{\partial B_1}(y)= \delta_{x_0}(y)$. From the Poisson formula, we have $u(y)= P(x_0, y)$.
}
\qed

\section{Gradient estimate in geodesic balls}

The following lemma is well-known. We include its proof here for completeness.
\begin{lemma}\label{lem Delta of Q}
{Assume $u: \Omega\rightarrow \mathbb{R}^+$ is a harmonic function, where $\Omega\subset M^n$ with $Rc\geq 0$, then for $Q= |\nabla \ln u|^2$, we have
\begin{align}
\Delta Q\geq \frac{2}{n- 1}Q^2+ \frac{4- 2n}{n- 1}\langle \nabla Q , \nabla v \rangle+ \frac{n}{2(n- 1)}|\nabla Q|^2\cdot Q^{- 1}. \nonumber 
\end{align}
}
\end{lemma}

\pf
{One can set $v= \ln u, Q= |\nabla v|^2$. Using the Bochner formula, we compute
\begin{align}
\Delta Q= 2v_{ij}^2+ 2Rc(\nabla v, \nabla v)+ 2\langle \nabla v, \nabla \Delta v\rangle \geq 2v_{ij}^2- 2\langle \nabla v, \nabla Q\rangle \label{Delta Q 1-new}
\end{align}
where we use the fact $\Delta v= \nabla (\frac{\nabla u}{u})= -|\nabla v|^2= -Q$.

Firstly observed by Yau \cite{Yau}, there is the following estimate for $\displaystyle \sum_{ij}v_{ij}^2$,
\begin{align}
v_{ij}^2&\geq v_{11}^2+ 2\sum_{\alpha= 2}^n v_{1\alpha}^2+ \sum_{\alpha= 2}^n v_{\alpha \alpha}^2 \geq v_{11}^2+ 2\sum_{\alpha= 2}^n v_{1\alpha}^2+ \frac{\big(\sum_{\alpha= 2}^n v_{\alpha \alpha}\big)^2}{n- 1} \nonumber \\
&= v_{11}^2+ 2\sum_{\alpha= 2}^n v_{1\alpha}^2+ \frac{\big(\Delta v- v_{11}\big)^2}{n- 1}= v_{11}^2+ 2\sum_{\alpha= 2}^n v_{1\alpha}^2+ \frac{\big(Q+ v_{11}\big)^2}{n- 1} \nonumber \\
&\geq \frac{n}{n- 1}\sum_{j= 1}^n v_{1j}^2+ \frac{1}{n- 1}Q^2+ \frac{2}{n- 1}Q\cdot v_{11} \label{square of vij-new}
\end{align}

Choosing an orthonormal frame $\{e_1, \cdots, e_n\}$ at a point such that $e_1= \frac{\nabla v}{|\nabla v|}$, 
\begin{align}
|\nabla Q|^2= \big|\nabla |\nabla v|^2\big|^2= 4\sum_{j= 1}^n \big(\sum_{i= 1}^n v_i v_{ij}\big)^2= 4v_1^2 \sum_{j= 1}^n v_{1j}^2= 4|\nabla v|^2 \sum_{j= 1}^n v_{1j}^2= 4Q \cdot \sum_{j= 1}^n v_{1j}^2\nonumber
\end{align}
which implies 
\begin{align}
\sum_{j= 1}^n v_{1j}^2= \frac{|\nabla Q|^2}{4Q} \label{square of v1j-new}
\end{align}

Combining the identity $\langle \nabla |\nabla v|^2, \frac{\nabla v}{|\nabla v|} \rangle=  \nabla_{e_1} |\nabla v|^2= 2v_1 v_{11}$ with $v_1= |\nabla v|$, we get 
\begin{align}
v_{11}= \frac{1}{2v_1|\nabla v|}\langle \nabla |\nabla v|^2 , \nabla v \rangle= \frac{1}{2Q}\langle \nabla Q , \nabla v \rangle \label{v11-new}
\end{align}

From (\ref{square of vij-new}), (\ref{square of v1j-new}) and (\ref{v11-new}), 
\begin{align}
v_{ij}^2\geq \frac{n}{4(n- 1)}|\nabla Q|^2\cdot Q^{- 1}+ \frac{1}{n- 1}Q^2+ \frac{1}{n- 1}\langle \nabla Q , \nabla v \rangle \nonumber 
\end{align}
combining the above estimate with (\ref{Delta Q 1-new}), we obtain
\begin{align}
\Delta Q\geq \frac{2}{n- 1}Q^2+ \frac{4- 2n}{n- 1}\langle \nabla Q , \nabla v \rangle+ \frac{n}{2(n- 1)}|\nabla Q|^2\cdot Q^{- 1} \nonumber 
\end{align}
}
\qed

From the above lemma and the maximum principle, we get the following corollary.
\begin{cor}\label{cor max of log gradient is on boundary}
{For a positive harmonic function $u$ on $B_1(p)$ with $Rc\geq 0$, we have 
\begin{align}
\sup\limits_{B_t(0)} \frac{|\nabla u|}{u}= \sup\limits_{\partial B_t(0)} \frac{|\nabla u|}{u} , \quad \quad \quad \quad \forall t\in (0, 1)\nonumber 
\end{align}
}
\end{cor}\qed

\begin{lemma}\label{lem Cheng-Yau gradient est by cut-off}
{For $B_1(p)\subseteq M^n$ with $Rc\geq 0$, assume $u: B_{1}(p)\rightarrow \mathbb{R}^+$ is a positive harmonic function, then for any positive function $h\in C^\infty(B_{1-\epsilon}(p))$ satisfying $\lim\limits_{r\rightarrow 1- \epsilon}\inf\limits_{x\in \partial B_r(p)}h(x)= \infty$, where $\epsilon> 0$ is a constant, we have
\begin{align}
\sup_{B_s(p)}\frac{|\nabla u|}{u}\leq \sup_{\partial B_s(p)} \sqrt{h}\cdot \sup_{y\in B_{1- \epsilon}(p)}T(h)(y), \quad \quad \quad \quad \forall s< 1-\epsilon \nonumber 
\end{align}
where  $T(h)= \frac{n- 2}{2}\frac{|\nabla h|}{h^{\frac{3}{2}}}+ \frac{1}{2}\sqrt{\Big\{(n^2- 5n+ 4)\frac{|\nabla h|^2}{h^3}+ 2(n- 1)\frac{\Delta h}{h^2}\Big\}^+}$.
}
\end{lemma}

\pf
{Let $\phi= \frac{1}{h}$ on $B_{1- \epsilon}(p)$ and $\phi= 0$ on $M^n- B_{1- \epsilon}(p)$, then $\phi$ is a positive cut-off function on $M^n$. Define $P= \phi Q= \phi |\nabla \ln u|^2$, from Lemma \ref{lem Delta of Q},
\begin{align}
&\quad \Delta P= Q\Delta \phi+ 2\langle \nabla \phi, \nabla Q\rangle+ \phi \Delta Q \nonumber \\
&\geq  \frac{\Delta \phi}{\phi}P+ 2\phi^{-1}\langle \nabla \phi, \nabla P\rangle- 2\phi^{-2}|\nabla \phi|^2\cdot P+ \frac{n}{2(n- 1)}|\nabla P|^2\cdot P^{-1} \nonumber \\
&+ \frac{n}{2(n- 1)}|\nabla \phi|^2\phi^{-2}\cdot P \nonumber -\frac{n}{n- 1}\phi^{-1}\langle \nabla \phi, \nabla P \rangle- \frac{2(n- 2)}{n- 1}\langle \nabla v, \nabla P \rangle \nonumber \\
&+ \frac{2(n- 2)}{n- 1}\langle \nabla v, \nabla \phi\rangle Q+ \frac{2}{n- 1}\phi^{-1}P^2 \nonumber 
\end{align}
Using the inequality $\big|\langle \nabla v, \nabla \phi\rangle\big|\leq |\nabla \phi|\cdot |\nabla v|=|\nabla \phi|\cdot \phi^{- \frac{1}{2}}P^{\frac{1}{2}}$, we get 
\begin{align}
\Delta P&\geq \frac{\Delta \phi}{\phi}P+ \frac{n -2}{n- 1}\phi^{-1}\langle \nabla \phi, \nabla P\rangle+ \frac{4- 3n}{2(n- 1)}\phi^{-2}|\nabla \phi|^2\cdot P+ \frac{n}{2(n- 1)}|\nabla P|^2\cdot P^{-1} \nonumber \\
&\quad - \frac{2(n- 2)}{n- 1}\langle \nabla v, \nabla P \rangle- \frac{2(n- 2)}{n- 1}|\nabla \phi|\phi^{-\frac{3}{2}}P^{\frac{3}{2}}+ \frac{2}{n- 1}\phi^{-1}P^2 \label{one side ineq}
\end{align}

Let $\mathcal{L}_1(P)\vcentcolon = \Delta P+ \frac{2(n- 2)}{n- 1}\langle \nabla v, \nabla P \rangle- \frac{n -2}{n- 1}\phi^{-1}\langle \nabla \phi, \nabla P\rangle-  \frac{n}{2(n- 1)}|\nabla P|^2\cdot P^{-1}$, then (\ref{one side ineq}) is equivalent to
\begin{align}
\mathcal{L}_1(P)&\geq \Big(\frac{\Delta \phi}{\phi}+ \frac{4- 3n}{2(n- 1)}\phi^{-2}|\nabla \phi|^2\Big)\cdot P  - \frac{2(n- 2)}{n- 1}|\nabla \phi|\phi^{-\frac{3}{2}}P^{\frac{3}{2}}+ \frac{2}{n- 1}\phi^{-1}P^2 \nonumber \\
&\geq \frac{1}{(n- 1)\phi}\Big\{(n- 1)P\Delta \phi+ \frac{4- 3n}{2}\phi^{-1}|\nabla \phi|^2\cdot P- 2(n- 2)|\nabla \phi|\phi^{-\frac{1}{2}}P^{\frac{3}{2}}+ 2P^2\Big\} \nonumber 
\end{align}

Assume $P(x_0)= \max\limits_{B_1(p)}P(x)$, then $x_0\in B_{1- \epsilon}(p)$ and $\mathcal{L}_1(P)(x_0)\leq 0$, we get 
\begin{align}
0\geq 2P(x_0)- 2(n- 2)\frac{|\nabla \phi|}{\sqrt{\phi}}P^{\frac{1}{2}}(x_0)+ \Big((n- 1)\Delta \phi+ \frac{4- 3n}{2}\phi^{-1}|\nabla \phi|^2\Big)(x_0) \nonumber 
\end{align}
which implies
\begin{align}
P(x_0)^{\frac{1}{2}}\leq \frac{n- 2}{2}\frac{|\nabla \phi|}{\sqrt{\phi}}(x_0)+ \frac{1}{2}\sqrt{(n^2- n)\frac{|\nabla \phi|^2}{\phi}- 2(n- 1)\Delta\phi}(x_0) \nonumber 
\end{align}

Putting $\phi= h^{-1}$ into the above, we have
\begin{align}
P(x_0)^{\frac{1}{2}}\leq \frac{n- 2}{2}\frac{|\nabla h|}{h^{\frac{3}{2}}}(x_0)+ \frac{1}{2}\sqrt{(n^2- 5n+ 4)\frac{|\nabla h|^2}{h^3}+ 2(n- 1)\frac{\Delta h}{h^2}}(x_0)= T(h)(x_0) \nonumber 
\end{align}

Then we get
\begin{align}
\sup_{B_{1- \epsilon}(p)} h(x)^{-\frac{1}{2}}\frac{|\nabla u|}{u}(x)= P(x_0)^{\frac{1}{2}}\leq T(h)(x_0)\leq \sup_{x\in B_{1- \epsilon}(p)}T(h)(x) \nonumber 
\end{align}
which implies $\frac{|\nabla u|}{u}(x)\leq \sqrt{h(x)}\sup\limits_{y\in B_{1- \epsilon}(p)}T(h)(y)$ for any $x\in B_{1- \epsilon}(p)$. Taking the supermum of $x$ on $\partial B_s(p)$, we get 
\begin{align}
\sup_{\partial B_s(p)}\frac{|\nabla u|}{u}(x)\leq \sup_{\partial B_s(p)}\sqrt{h(x)}\cdot \sup\limits_{y\in B_{1- \epsilon}(p)}T(h)(y), \quad \quad \quad \quad \forall s< 1- \epsilon \nonumber 
\end{align}
The conclusion follows from the above and Corollary \ref{cor max of log gradient is on boundary}.
}
\qed

\begin{theorem}\nonumber 
{Assume $u: B_{1}(p)\rightarrow \mathbb{R}^+$ is a harmonic function, where $M^2$ is a Riemannian manifold with $Rc\geq 0$, then $\sup\limits_{x\in B_s(p)}\frac{|\nabla u|}{u}\leq \frac{1}{1- s}+ \frac{1}{1+ s}$ for any $s\in [0, 1)$.
}
\end{theorem}

\pf
{When $n= 2$, choosing $h(x)= \frac{4}{(1- \epsilon- \rho(x)^2)^2}$ in Lemma \ref{lem Cheng-Yau gradient est by cut-off}, we get
\begin{align}
\sup_{B_s(p)}\frac{|\nabla u|}{u}\leq \frac{1}{2}\sup_{\partial B_s(p)} \sqrt{h}\cdot \sup_{y\in B_{1- \epsilon}(p)}\sqrt{2\frac{\Delta h}{h^2}- 2\frac{|\nabla h|^2}{h^3}} , \quad \quad \quad \quad \forall s< 1-\epsilon\nonumber 
\end{align}
By the Laplace comparison theorem, we have $\Delta \rho\leq \frac{1}{\rho}$, which yields $\sqrt{2\frac{\Delta h}{h^2}- 2\frac{|\nabla h|^2}{h^3}}\leq 2$. Hence we get
\begin{align}
\sup_{B_s(p)}\frac{|\nabla u|}{u}\leq \sup_{\partial B_s(p)} \sqrt{h}= \frac{2}{(1- \epsilon- s^2)}  , \quad \quad \quad \quad \forall s< 1-\epsilon\nonumber 
\end{align}
Letting $\epsilon\rightarrow 0$ in the above, the conclusion follows.
}
\qed

The following corollary follows from the above theorem directly.
\begin{cor}\label{cor growth of harm func on surface}
{Assume $u: B_{1}(p)\rightarrow \mathbb{R}^+$ is a harmonic function, where $M^2$ is a Riemannian manifold with $Rc\geq 0$, then $\displaystyle u(x)\leq \frac{1+ d(p, x)}{1- d(p, x)}u(p)$.
}
\end{cor}\qed

\begin{remark}\label{rem rigidity for 2-dim}
{We do not know the rigidity part of the above Theorem. In other words, if $\sup\limits_{x\in B_s(p)}\frac{|\nabla u|}{u}=\frac{1}{1- s}+ \frac{1}{1+ s}$ for some $s\in [0, 1)$, is $B_1(p)$ isometric to $B_1(0)$?
}
\end{remark}

\begin{question}\label{ques sharp est of Poisson kernel}
{Assume $u: B_{1}(p)\rightarrow \mathbb{R}^+$ is a positive harmonic function, where $B_1(p)\subseteq M^n$ with $Rc\geq 0$ and $n\geq 3$, for any $s\in [0, 1)$, do we have
\begin{align}
\sup\limits_{B_s(p)}\frac{|\nabla u|}{u}\leq \frac{n- 1}{1- s}+ \frac{1}{1+ s} ? \nonumber 
\end{align}
}
\end{question}


\section*{Part II. Polynomial growth harmonic functions on manifolds}

In part II of this paper, we always assume that $M^n$ is a complete non-compact Riemannian manifold with $Rc\geq 0$ and maximal volume growth, where $n\geq 3$, unless otherwise mentioned.

\section{Frequency and existence of harmonic function}

From \cite{Var}, there exists a unique, minimal, positive Green function on $M^n$, denoted as $G(p, x)$ where $p$ is a fixed point on manifold. We define $\displaystyle b(x)= [n(n- 2)\omega_n\cdot G(p, x)]^{\frac{1}{2- n}}$ and use $\displaystyle \rho(x)= d(p, x)$ unless otherwise mentioned.  From the behavior of Green function $G(p, x)$ near singular point $p$, we get 
\begin{align}
\lim_{\rho(x)\rightarrow 0}\frac{b(x)}{\rho(x)}= 1 \quad \quad and \quad \quad \lim_{\rho(x)\rightarrow 0} |\nabla b|= 1. \label{est of gradient of b at 0}
\end{align}

The following definition of frequency function was firstly introduced in \cite{CM-JDG-97}. Assume $u(x)$ is a harmonic function defined on $\{b(x)\leq r\}$, define:
\begin{align}
I_{u}(r)= r^{1- n}\int_{b(x)= r} u^2|\nabla b| dx \quad \quad and \quad \quad D_{u}(r)= r^{2- n}\int_{b(x)\leq r} |\nabla u|^2 dx  .\label{def of D(r)} 
\end{align}
\textbf{The frequency function of harmonic function $u$} is defined by $\mathscr{F}_{u}(r)= \frac{D_{u}(r)}{I_{u}(r)}$. 

In this paper, $I(r)$, $D(r)$ and $\mathscr{F}(r)$ are defined as in above with respect to harmonic function $u$ (which may be defined locally on suitable regions of $M^n$). From the definition of $\mathscr{F}_u(r)$, we can get that $\mathscr{F}_u(0)= 0$ if $u(p)\neq 0$.

The following lemma is a generalized version of the Rellich-Necas identity, which was originally discovered by Payne and Weinberger \cite{PW} (also see \cite{CM-JDG-97}).

\begin{lemma}\label{lem mysterious formula}
{For harmonic function $u$ on $B= \{x: b(x)\leq r\}\subseteq M^n$, we have
\begin{align}
r\int_{b= r} |\nabla u|^2\cdot |\nabla b|= 2r\int_{b= r} \Big|\frac{\partial u}{\partial \vec{n}}\Big|^2\cdot |\nabla b|+ \int_{b\leq r} \frac{1}{2}\Delta (b^2)\cdot |\nabla u|^2- \nabla^2(b^2)(\nabla u, \nabla u) ,\nonumber 
\end{align}
where $\vec{n}$ is the unit normal of $\partial B$ pointing inward of $B$.
}
\end{lemma}\qed

\begin{lemma}\label{lem lower bound of freq's deri}
{When $I(r)\neq 0$, we have
\begin{align}
\big(\ln\mathscr{F}(r)\big)'\geq -\frac{\int_{b\leq r} \Big\{n\big(1- |\nabla b|^2\big)+ \big|2g- \nabla^2(b^2)\big|\Big\}\cdot |\nabla u|^2}{r\int_{b\leq r}|\nabla u|^2}- 2\frac{\int_{b= r} |\frac{\partial u}{\partial \vec{n}}|^2\big||\nabla b|- |\nabla b|^{-1}\big|}{\int_{b\leq r} |\nabla u|^2}. \nonumber
\end{align}
}
\end{lemma}

\pf
{From the definition of $\mathscr{F}$, we have 
\begin{align}
\mathscr{F}'(r)= \frac{D'(r)}{I(r)}- \frac{D(r)}{I^2(r)}I'(r)= \mathscr{F}(r)\Big\{\frac{D'(r)}{D(r)}- \frac{I'(r)}{I(r)}\Big\}. \nonumber 
\end{align}

From Lemma \ref{lem derivative of I(r)} and Cauchy-Schwartz inequality, noting that $\int_{b\leq r}|\nabla u|^2= \int_{b= r} u\frac{\partial u}{\partial \vec{n}}$, we have 
\begin{align}
&\quad \frac{D'(r)}{D(r)}- \frac{I'(r)}{I(r)}= \frac{2- n}{r}+ \frac{\int_{b= r}|\nabla u|^2 |\nabla b|^{-1}}{\int_{b\leq r}|\nabla u|^2}- 2\frac{\int_{b\leq r} |\nabla u|^2}{\int_{b= r} u^2|\nabla b|} \nonumber \\
&\geq \frac{2- n}{r}+ \frac{\int_{b= r}|\nabla u|^2 |\nabla b|^{-1}}{\int_{b\leq r}|\nabla u|^2}- 2\frac{\int_{b= r} |\frac{\partial u}{\partial \vec{n}}|^2|\nabla b|^{-1}}{\int_{b\leq r} |\nabla u|^2} \nonumber\\
&= \frac{2- n}{r}+ \frac{\int_{b= r}|\nabla u|^2 |\nabla b|^{-1}- 2|\frac{\partial u}{\partial \vec{n}}|^2|\nabla b|}{\int_{b\leq r}|\nabla u|^2}  + 2\Big\{\frac{\int_{b= r} |\frac{\partial u}{\partial \vec{n}}|^2|\nabla b|}{\int_{b\leq r} |\nabla u|^2}- \frac{\int_{b= r} |\frac{\partial u}{\partial \vec{n}}|^2|\nabla b|^{-1}}{\int_{b\leq r} |\nabla u|^2}\Big\}. \nonumber 
\end{align}

Note $|\nabla b|\leq 1$ and Lemma \ref{lem mysterious formula}, we have
\begin{align}
&\quad \frac{D'(r)}{D(r)}- \frac{I'(r)}{I(r)}\geq \frac{2- n}{r}+ \frac{\int_{b= r}|\nabla u|^2 |\nabla b|- 2|\frac{\partial u}{\partial \vec{n}}|^2|\nabla b|}{\int_{b\leq r}|\nabla u|^2} + 2\frac{\int_{b= r} |\frac{\partial u}{\partial \vec{n}}|^2\big(|\nabla b|- |\nabla b|^{-1}\big)}{\int_{b\leq r} |\nabla u|^2} \nonumber \\
&= \frac{2- n}{r}+ \frac{1}{r}\cdot \frac{\int_{b\leq r} \Big[\frac{1}{2}\Delta(b^2)\cdot |\nabla u|^2- \nabla^2(b^2)(\nabla u, \nabla u)\Big]}{\int_{b\leq r}|\nabla u|^2}+ 2\frac{\int_{b= r} |\frac{\partial u}{\partial \vec{n}}|^2\big(|\nabla b|- |\nabla b|^{-1}\big)}{\int_{b\leq r} |\nabla u|^2} \nonumber \\
&= \frac{\int_{b\leq r} n\Big[|\nabla b|^2- 1\Big]|\nabla u|^2+ \Big[\big(2g- \nabla^2(b^2)\big)(\nabla u, \nabla u)\Big]}{r\int_{b\leq r}|\nabla u|^2}+ 2\frac{\int_{b= r} |\frac{\partial u}{\partial \vec{n}}|^2\big(|\nabla b|- |\nabla b|^{-1}\big)}{\int_{b\leq r} |\nabla u|^2} \nonumber\\
&\geq -\frac{\int_{b\leq r} \Big\{n\big(1- |\nabla b|^2\big)+ \big|2g- \nabla^2(b^2)\big|\Big\}\cdot |\nabla u|^2}{r\int_{b\leq r}|\nabla u|^2}- 2\frac{\int_{b= r} |\frac{\partial u}{\partial \vec{n}}|^2\big||\nabla b|- |\nabla b|^{-1}\big|}{\int_{b\leq r} |\nabla u|^2} ,\nonumber
\end{align}
in the last equation above we use $\Delta (b^2)= 2n|\nabla b|^2$.
}
\qed

\begin{lemma}\label{lem set inclusion between b and rho}
{For any $\epsilon_0\in (0, 1)$, there is $r_0> 0$ such that if $r\geq \mathrm{V}_M^{\frac{1}{n- 2}}(1- \epsilon_0)^{\frac{1}{2- n}}r_0$, then
\begin{align}
\Big\{x: \rho(x)\leq \mathrm{V}_M^{\frac{1}{2- n}}(1- \epsilon_0)^{\frac{1}{2- n}}r\Big\}\subseteq \{x: b(x)\leq r\}\subseteq \Big\{x: \rho(x)\leq \mathrm{V}_M^{\frac{1}{2- n}}(1+ \epsilon_0)^{\frac{1}{n- 2}}r\Big\} .\nonumber 
\end{align} 
}
\end{lemma}

\pf
{From the proof of \cite[Theorem $1.1$]{LTW} (also see \cite{CM}), for any $\delta \in (0, \frac{1}{2}]$, we have 
\begin{align}
\big(\mathrm{V}_M\big)^{\frac{1}{n- 2}}\Big(1+ \tau(n, \delta)\Big)^{\frac{1}{2- n}}\rho(x) \leq b(x)\leq \big(\mathrm{V}_M\big)^{\frac{1}{n- 2}}\Big(1- \tau(n, \delta)\Big)^{\frac{1}{2- n}}\rho(x), \nonumber 
\end{align}
where $\tau(n, \delta)(x)= C(n)\big[\delta+ (\theta_p(\delta \rho(x))- \theta)^{\frac{1}{n- 1}}\big]$, and 
\begin{align}
\theta_p(r)= \frac{\mathrm{Vol}(\partial B_r(p))}{r^{n- 1}}, \quad \quad \quad \theta= \lim\limits_{r\rightarrow \infty}\theta_p(r).\nonumber 
\end{align}

For any$\displaystyle\epsilon_0\in (0, 1)$, we can firstly find $\displaystyle\delta= \frac{\epsilon_0}{2C(n)}$, then
\begin{align}
\tau(n, \delta)= C(n)\delta+ C(n)(\theta_p(\delta \rho(x))- \theta)^{\frac{1}{n- 1}}= \frac{\epsilon_0}{2}+ C(n)\Big\{\theta_p(\frac{\epsilon_0}{2C(n)} \rho(x))- \theta\Big\}^{\frac{1}{n- 1}}. \nonumber 
\end{align}

Now choose suitable $\displaystyle r_0> 0$, then it yields $\displaystyle \tau(n, \delta)(x)\leq \epsilon_0$. If $\displaystyle \rho(x)\geq r_0$, we have 
\begin{align}
\big(\mathrm{V}_M\big)^{\frac{1}{n- 2}}\Big(1+ \epsilon_0\Big)^{\frac{1}{2- n}}\rho(x) \leq b(x)\leq \big(\mathrm{V}_M\big)^{\frac{1}{n- 2}}\Big(1- \epsilon_0\Big)^{\frac{1}{2- n}}\rho(x). \label{lower and upper bound of b} 
\end{align}

From \cite{Colding}, we know that $\displaystyle |\nabla b|\leq 1$. If $\displaystyle \rho(x)\leq r_0$, we get $\displaystyle b(x)\leq \rho(x)\leq r_0\leq r$. If  $\displaystyle \rho(x)\in [r_0, \mathrm{V}_M^{\frac{1}{2- n}}(1- \epsilon_0)^{\frac{1}{2- n}}r]$, then from (\ref{lower and upper bound of b}), 
\begin{align}
b(x)\leq \mathrm{V}_M^{\frac{1}{n- 2}}(1- \epsilon_0)^{\frac{1}{2- n}}\rho(x)\leq r, \nonumber 
\end{align}
the first inclusion is proved.

If $\displaystyle b(x)\leq r$, and  $\displaystyle \rho(x)> (1+ \epsilon_0)^{\frac{1}{n- 2}}\mathrm{V}_M^{\frac{1}{2- n}}r$, then  $\displaystyle \rho(x)> r_0$. From (\ref{lower and upper bound of b}),
\begin{align}
b(x)\geq \mathrm{V}_M^{\frac{1}{n- 2}}(1+ \epsilon_0)^{\frac{1}{2- n}}\rho(x)> r. \nonumber 
\end{align}
This is the contradiction, and the second inclusion follows.
}
\qed

\begin{lemma}\label{lem varia of mean-value ineq}
{There is  $\displaystyle R= R(M^n, n)> 0$, such that for any  $\displaystyle \tau> 1, r\geq R$, we have 
\begin{align}
\sup_{b(x)\leq r} |\nabla u|^2(x)\leq \frac{C(n, \tau)}{\mathrm{V}_M^{\frac{2}{2- n}}}r^{-2} D(\tau r) \quad \quad and \quad \quad \sup\limits_{b\leq r}|u|^2\leq \frac{C(n, \tau)}{\mathrm{V}_M^{\frac{2}{2- n}}}D(\tau r), \nonumber 
\end{align}
where  $\displaystyle u$ is harmonic with  $\displaystyle u(p)= 0$.
}
\end{lemma}

\pf
{From Lemma \ref{lem set inclusion between b and rho}, choose  $\displaystyle \epsilon_0= \frac{1}{\sqrt{2}}$, there is  $\displaystyle r_0> 0$ such that if  $\displaystyle r\geq \mathrm{V}_M^{\frac{1}{n- 2}}(1- \epsilon_0)^{\frac{1}{2- n}}r_0$, then
\begin{align}
\Big\{\rho\leq \mathrm{V}_M^{\frac{1}{2- n}}(1- \epsilon_0)^{\frac{1}{2- n}}r\Big\}\subseteq \{b\leq r\}\subseteq \Big\{\rho\leq \mathrm{V}_M^{\frac{1}{2- n}}(1+ \epsilon_0)^{\frac{1}{n- 2}}r\Big\} .\label{set inclusion} 
\end{align} 

From \cite[Theorem $1.2$]{LS} and (\ref{set inclusion}), let  $\displaystyle \tau_1= 2^{\frac{1}{n- 2}}\tau> 0$, 
\begin{align}
\sup_{b(x)\leq r} |\nabla u|^2(x)&\leq \sup_{B(\mathrm{V}_M^{\frac{1}{2- n}}(1+ \epsilon_0)^{\frac{1}{n- 2}}r)} |\nabla u|^2 \nonumber \\
&\leq \frac{C(n, \tau_1)}{\mathrm{Vol}(B(\tau_1\mathrm{V}_M^{\frac{1}{2- n}}(1+ \epsilon_0)^{\frac{1}{n- 2}}r))} \int_{B(\tau_1\mathrm{V}_M^{\frac{1}{2- n}}(1+ \epsilon_0)^{\frac{1}{n- 2}}r)} |\nabla u|^2  \nonumber \\
&\leq \frac{C(n, \tau)}{\mathrm{V}_M\cdot \Big(\tau_1\mathrm{V}_M^{\frac{1}{2- n}}(1+ \epsilon_0)^{\frac{1}{n- 2}}r\Big)^n} \int_{b\leq \tau_1\big(1- \epsilon_0^2\big)^{\frac{1}{n- 2}}r} |\nabla u|^2 \nonumber \\
&\leq \frac{C(n, \tau)}{\mathrm{V}_M^{\frac{2}{2- n}}}r^{-2} D(\tau r).  \label{4.4.3}
\end{align}

Integrating (\ref{4.4.3}) along geodesics starting at  $\displaystyle p$ and using  $\displaystyle u(p)= 0$, we obtain  $\displaystyle \sup\limits_{b(x)\leq r} |u|^2(x)\leq \frac{C(n, \tau)}{\mathrm{V}_M^{\frac{2}{2- n}}} D(\tau r)$.
}
\qed

\begin{lemma}\label{lem poly growth by frequency}
{If  $\displaystyle u(x)$ is harmonic on  $\displaystyle B_{2r}(p)\subset M^n$ with  $\displaystyle \max\limits_{s\leq 2r}\mathscr{F}_u(s)\leq d$ and  $\displaystyle u(p)= 0$, where  $\displaystyle r\geq \max\{1, R(M^n, p, n)\}$, then  $\displaystyle \sup\limits_{b\leq r}|u|\leq C(n, \mathrm{V}_M, d)r^{d}I_u(1)^{\frac{1}{2}}$.
}
\end{lemma}

\pf
{From Lemma \ref{lem varia of mean-value ineq}, we have
\begin{align}
\sup_{b\leq r}|u|^2\leq \frac{C(n)}{\mathrm{V}_M^{\frac{2}{2- n}}}D(2r), \quad \quad \quad \quad \forall r\geq R(M^n, p, n). \label{upper bound of u} 
\end{align}

By Lemma \ref{lem derivative of I(r)}, for any $s\in [1, 2r]$, we have 
\begin{align}
I(s)= \exp\Big(2\int_1^s \frac{\mathscr{F}(t)}{t}dt\Big)I(1)\leq I(1) s^{2d}.\label{upper bound of I}
\end{align}
From (\ref{upper bound of I}), we have 
\begin{align}
D(s)= I(s)\mathscr{F}(s)\leq I(1)d\cdot s^{2d}, \quad \quad \quad \quad \forall \ s\in [1, 2r]. \label{upper bound of D}
\end{align}

By (\ref{upper bound of D}) and (\ref{upper bound of u}), we obtain
\begin{align}
\sup_{b\leq r}|u|^2&\leq \frac{C(n)D(2r)}{\mathrm{V}_M^{\frac{2}{2- n}}}\leq C(n, \mathrm{V}_M)I(1)dr^{2d}= C(n, \mathrm{V}_M, d)r^{2d}I(1) .\nonumber 
\end{align}
}
\qed

\begin{lemma}\label{lem criteria for exsitence of harm func of poly growth}
{If there exists a sequence of non-zero harmonic functions $u_i(x)$ defined on $B_{r_i}(p)\subseteq M^n$, where $\lim\limits_{i\rightarrow \infty}r_i= \infty$, satisfying $u_i(p)= 0$ and $\max\limits_{s\leq r_i}\mathscr{F}_{u_i}(s)\leq d$ for some $d> 0$, then there exists a non-constant polynomial growth harmonic function $u$ with degree $\leq d$.
}
\end{lemma}

\pf
{Consider $I_{u_i}(1)$, if $I_{u_i}(1)= 0$, from the maximum principle, we know that $u_i\big|_{B_1(p)}= 0$. Applying unique continuation theorem on $u_i$, we get $u_i\big|_{B_{r_i}(p)}= 0$, which is the contradiction. Hence $I_{u_i}(1)\neq 0$, we can define $\tilde{u}_i(x)= \frac{u_i(x)}{I_{u_i}(1)^{\frac{1}{2}}}$.

From Lemma \ref{lem derivative of I(r)}, we have $n\omega_nu_i(p)^2= \lim\limits_{r\rightarrow 0}I_{u_i}(r)\leq I_{u_i}(1)$, which implies
\begin{align}
|u_i(p)|\leq \Big(\frac{I_{u_i}(1)}{n\omega_n}\Big)^{\frac{1}{2}}. \label{upper bound of u_i at p}
\end{align}

From (\ref{upper bound of u_i at p}), we have 
\begin{align}
|\tilde{u}_i(p)|\leq \Big(\frac{1}{n\omega_n}\Big)^{\frac{1}{2}},\quad \quad \quad \quad I_{\tilde{u}_i}(1)= 1. \label{ui is not constant} 
\end{align}

From Lemma \ref{lem poly growth by frequency} and the assumption on $u_i$, using (\ref{upper bound of u_i at p}), we know that 
\begin{align}
\sup_{B_{\frac{r_i}{2}}(p)} |u_i|\leq |u_i(p)|+ C(n, d)r_i^dI_{u_i}(1)^{\frac{1}{2}}\leq C(n, d)r_i^dI_{u_i}(1)^{\frac{1}{2}} ,\nonumber 
\end{align}
which implies $\sup\limits_{B_{\frac{r_i}{2}}(p)} |\tilde{u}_i|\leq C(n, d)r_i^d$. From compactness theorem of harmonic functions and the above growth estimate for $u_i$, after taking suitable subsequence, $\tilde{u}_i$ converges to a polynomial growth harmonic function $u(x)$ with degree $\leq d$ on $M^n$. From (\ref{ui is not constant}), we know that $u$ satisfies $I_{u}(1)= 1$, hence $u(x)$ is not constant by $u_i(p)= 0$. The conclusion is proved.
}
\qed

\begin{prop}\label{prop local bound of freq}
{For $\gamma\geq 1$, there is $R_2= R_2(M^n, p, \delta, \gamma)$, such that for $R\geq R_2$ and harmonic function $u$ on $B_R\subseteq M^n$, if $I_{u}(r)\leq \gamma \cdot I_{u}(\frac{r}{2})$ for any $r\in [2^{-1} R, 2 R]$, then $\sup\limits_{s\in [\frac{1}{2}R, R]}\mathscr{F}(s)\leq 14\ln\gamma$.
}
\end{prop}

\pf
{\textbf{Step (1)}.From Lemma \ref{lem varia of mean-value ineq}, we obtain $\sup\limits_{b(x)\leq \frac{r}{2}} |u|^2(x)\leq \frac{C(n)}{\mathrm{V}_M^{\frac{2}{2- n}}} D(\frac{5}{8} r)$. Hence
\begin{align}
I(\frac{r}{2})&= \big(\frac{r}{2}\big)^{1- n}\int_{b= \frac{r}{2}} u^2|\nabla b|\leq \frac{C(n)}{\mathrm{V}_M^{\frac{2}{2- n}}} D(\frac{5}{8} r)\cdot \big(\frac{r}{2}\big)^{1- n} \int_{b= \frac{r}{2}} |\nabla b|\nonumber \\
&= \frac{C(n)}{\mathrm{V}_M^{\frac{2}{2- n}}} D(\frac{5}{8} r) I_1(\frac{r}{2}) = C(n)\mathrm{V}_M^{\frac{2}{n- 2}} D(\frac{5}{8} r) ,\label{4.4.4}
\end{align}
in the last equation we use Lemma \ref{lem derivative of I(r)}. Also from Lemma \ref{lem derivative of I(r)},
\begin{align}
\int_{\frac{14}{15}r}^{r} \frac{2D(s)}{s} ds= I(r)- I(\frac{14}{15}r)\leq I(r) ,\nonumber 
\end{align}
which implies $2\int_{\frac{14}{15}r}^{r} s^{n- 2}D(s) ds\leq r^{n- 1}I(r)$. Note $s^{n- 2}D(s)$ is nondecreasing in $s$ from the definition of $D(r)$, we get $\frac{1}{7}\Big(\frac{14}{15}r\Big)^{n- 1} D(\frac{14}{15}r)\leq r^{n- 1}I(r)$. Combining (\ref{4.4.4}), simplifying the above inequality yields 
\begin{align}
D(\frac{14}{15}r)\leq C(n) I(r)\leq C(n)\gamma I(\frac{r}{2})\leq C(n)\mathrm{V}_M^{\frac{2}{n- 2}}\gamma D(\frac{5}{8}r). \label{eq 0.1}
\end{align}

From Lemma \ref{lem varia of mean-value ineq}, let $\tau= \Big(\frac{16}{15}\Big)^{\frac{1}{3}}$, we have 
\begin{align}
\sup_{b(x)\leq \frac{7}{8}r} |\nabla u|^2(x)&\leq C(n)\mathrm{V}_M^{\frac{2}{n- 2}}r^{-2} D(\frac{14}{15} r)\leq C(n, \mathrm{V}_M)r^{-2}\gamma D(\frac{5}{8} r) \nonumber \\
&= C(n, \mathrm{V}_M)\gamma r^{-n}\int_{b\leq \frac{5}{8} r} |\nabla u|^2, \nonumber 
\end{align}
hence we get
\begin{equation}
\sup_{b(x)\leq s} |\nabla u|^2(x)\leq  C(n, \mathrm{V}_M)\gamma s^{-n}\int_{b\leq s}|\nabla u|^2 , \quad \quad \quad \quad \forall s\in [\frac{5}{8}r, \frac{7}{8}r]. \nonumber 
\end{equation}

\textbf{Step (2)}. From \cite{CC-Ann}, given any $\delta> 0$, there exists $R_1= R_1(M^n, p, \delta)> 0$ such that for $r\geq R_1$, we have
\begin{align}
\fint_{b(y)\leq r} \Big||\nabla b|^2- 1\Big| dy \leq \delta  \ ,  \quad \quad 
\fint_{b(y)\leq r} \Big|\nabla^2(b^2)- 2g\Big| dy \leq \delta. \label{almost cone}
\end{align}

In the rest of the proof, we assume that $R\geq 16R_1(M^n, p, \delta)$, where $\delta$ is to be determined later. Using (\ref{almost cone}), Lemma \ref{lem lower bound of freq's deri} and Step (1), for $s\in \Big[\frac{5}{8}r, \frac{7}{8}r\Big]$ we have 
\begin{align}
\Big(\ln \mathscr{F}(s)\Big)'&\geq -\frac{\sup\limits_{b\leq s}|\nabla u|^2}{\int_{b\leq s}|\nabla u|^2}\Big\{\frac{n}{s}\int_{b\leq s}\big||\nabla b|^2- 1\big|+ \big|2g- \nabla^2(b^2)\big|+ 2\int_{b= s}\big||\nabla b|- |\nabla b|^{-1}\big|\Big\}  \nonumber \\
&\geq -C(n, \mathrm{V}_M)\gamma \cdot \Big\{\frac{\int_{b\leq s}\big||\nabla b|^2- 1\big|+ \big|2g- \nabla^2(b^2)\big|}{s^{n+ 1}}+ \frac{\int_{b= s}\big||\nabla b|- |\nabla b|^{-1}\big|}{s^n}\Big\} \nonumber \\
&\geq -C(n, \mathrm{V}_M)\gamma \cdot\Big\{\frac{\fint_{b\leq s}\big||\nabla b|^2- 1\big|+ \big|2g- \nabla^2(b^2)\big|}{s}+ \frac{\int_{b= s}\big||\nabla b|- |\nabla b|^{-1}\big|}{r^n}\Big\} \nonumber \\
&\geq -C(n, \mathrm{V}_M)\gamma \cdot\Big\{\frac{\delta}{s}+ \frac{\int_{b= s}\big||\nabla b|- |\nabla b|^{-1}\big|}{r^n}\Big\}. \nonumber 
\end{align}

Hence take the integral of the above inequality, we obtain
\begin{align}
\int_{\frac{5}{8} r}^{\frac{7}{8} r} (\ln \mathscr{F}(s))' ds&\geq -C(n, \mathrm{V}_M) \gamma \delta - C(n, \mathrm{V}_M)\gamma r^{-n}\int_{\frac{5}{8}r\leq b\leq \frac{7}{8}r} \big||\nabla b|^2- 1\big| \nonumber \\
&\geq -C(n, \mathrm{V}_M) \gamma \delta - C(n, \mathrm{V}_M)\gamma \fint_{b\leq \frac{7}{8}r} \big||\nabla b|^2- 1\big|  \nonumber \\
&\geq -C(n, \mathrm{V}_M) \gamma\delta,  \nonumber
\end{align}
in the first inequality above we use the Co-area formula. Hence there exists $\delta_0= \delta_0(n, \mathrm{V}_M, \gamma)$, such that if $\delta\leq \delta_0$,
\begin{align}
\int_{\frac{5}{8} r}^{\frac{7}{8} r}\big(\ln \mathscr{F}\big)'(s) ds\geq -1 . \nonumber 
\end{align}

From Step (1) and Lemma \ref{lem derivative of I(r)}, 
\begin{align}
\int_{\frac{5}{8}r}^{\frac{7}{8}r} \frac{2\mathscr{F}(s)}{s} ds= \Big(\ln I(s)\Big)\Big|_{\frac{5}{8}r}^{\frac{7}{8}r}\leq \ln\frac{I(r)}{I(\frac{r}{2})}\leq \ln\gamma . \nonumber 
\end{align}
Hence there exists $s_i\in \big[\frac{3}{4}r, \frac{7}{8}r\big]$ such that $\mathscr{F}(s_i)\leq \frac{7}{2}\ln\gamma$. For this $s_i$, similar to the above argument, in fact we have $\int_{s}^{s_i}\big(\ln \mathscr{F}\big)'(t) dt\geq -1$ for any $s\in [\frac{5}{8}r, \frac{11}{16}r]$, which implies 
\begin{align}
\sup\limits_{s\in [\frac{5}{8}r, \frac{11}{16}r]}\mathscr{F}(s)\leq e\mathscr{F}(s_i)\leq 14\ln\gamma, \quad \quad \quad \quad \forall r\in [2^{-1} R, 2 R] \nonumber 
\end{align}
then we get that$\sup\limits_{s\in [\frac{1}{2}R, R]}\mathscr{F}(s)\leq e\mathscr{F}(s_i)\leq 14\ln\gamma$.  
}
\qed

Now we conclude this section with the proof of Theorem \ref{thm quc imply existence}.

\pf[of Theorem \ref{thm quc imply existence}]
{\textbf{Step (1)}. From \cite{CC-Ann}, we know any tangent cone at infinity of $M^n$ is a metric cone, choose one denoted as $C(X)$. Consider $\varphi_1(x)$ is the eigenfunction on $X$ with respect to eigenvalue $\lambda_1= \alpha_1(\alpha_1+ n- 2)$, $\int_{X} |\varphi_1|^2= 1$. Let $u_{\infty}= r^{\alpha_1} \varphi_1(x)$, then $u_{\infty}$ is harmonic on $C(X)$ and $u_{\infty}(p_{\infty})= 0$.

From \cite[Lemma $4.1$]{Xu}, there exist $R_i\rightarrow \infty$, $B(R_i)\subset M^n$, such that 
\begin{align}
\lim\limits_{i\rightarrow \infty} d_{GH} \big(B_i(1), B_{\infty}(1)\big)= 0 ,\nonumber 
\end{align}
where $B_i(1)\subset (M^n, R_i^{-2}g)$, $B_{\infty}(1)\subset C(X)$. And $u_i$ is harmonic on $B_p(R_i)= B_i(1)$ satisfying the following property:
\begin{align}
\lim_{i\rightarrow \infty} |u_i\circ \Psi_{\infty, i}- u_{\infty}|_{L^{\infty}\big(B_{\infty}(1)\big)}=0\ , \quad\quad\quad u_i(p)= 0 , \label{5.1.1}
\end{align}
where $\Psi_{\infty, i}: B_{\infty}(1)\rightarrow  B_i(1)$ is an $\epsilon_i$-Gromov-Hausdorff approximation, and $\lim\limits_{i\rightarrow \infty} \epsilon_i= 0$. From Proposition $3.4$ in \cite{Honda}, we have
\begin{align}
\lim_{i\rightarrow \infty} I_{u_i}(t)= I_{u_{\infty}}(t), \quad \quad \quad \quad \forall t\in (0, 1] .\label{3.5.8}
\end{align}

\textbf{Step (2)}. Let $d= \alpha_1+ 1$, we will prove that there exists $i_0= i_0> 0$ such that if $i\geq i_0$, then 
\begin{equation}
I_{u_i}(r)\leq 2^{2d} \cdot I_{u_i}(\frac{r}{2}) , \quad \quad \quad \quad \forall r\in [4^{-1}R_i, R_i] .\label{5.2.1}
\end{equation}

By contradiction. If the above statement does not hold, we can assume that there exists a subsequence $\{r_i\}$ with $r_i\in [4^{-1}R_i, R_i]$, such that
\begin{align}
I_{u_i}(r_i)> 2^{2d} \cdot I_{u_i}(\frac{r_i}{2}) . \label{5.2.2}
\end{align}

Note $R_i^{-1}r_i\in [4^{-1}, 1]$, without loss of generality, we can assume that there exists a subsequence of $\{i\}$, for simplicity also denoted as $\{i\}$ such that
\begin{align}
\lim_{i\rightarrow \infty}R_i^{-1}r_i= c_0\in [4^{-1}, 1] ,\label{5.2.3}
\end{align}
where $c_0$ is some constant. Taking the limit in (\ref{5.2.2}), from (\ref{5.2.3}) and (\ref{3.5.8}), note $b_{\infty}= \rho_{\infty}$ on $C(X)$, we have
\begin{align}
c_0^{1- n}\int_{\rho_\infty= c_0} |u_{\infty}|^2 d\nu_{\infty}\geq 2^{2d}\big(\frac{c_0}{2}\big)^{1- n} \int_{\rho_\infty= \frac{c_0}{2}} |u_{\infty}|^2 d\nu_{\infty}. \label{5.2.4}
\end{align}
From $u_{\infty}= r^{\alpha_1} \varphi_1(x)$ and $d> \alpha_1$, (\ref{5.2.4}) implies $\int_{X} |\varphi_1(x)|^2 dx= 0$, which is contradiction. 

From (\ref{5.2.1}) and Proposition \ref{prop local bound of freq}, we get that $\sup\limits_{s\in [4^{-1}R_i, 2^{-1}R_i]} \mathscr{F}_{u_i}(s)\leq 28d\ln 2$. From quantitative strong unique continuation, we have $\sup\limits_{s\in [0, 2^{-1}R_i]}\mathscr{F}_{u_i}(s)\leq C(n, d, \mathrm{V}_M)$. Now applying Lemma \ref{lem criteria for existence of harm func of poly growth}, we obtain the existence of non-constant polynomial growth harmonic function on $M^n$.
}
\qed

\section{The frequency of linear growth harmonic functions}

\begin{lemma}\label{lem derivative of I(r)}
{If $u$ is a harmonic function, we have $I'_{u}(r)= 2\frac{D_{u}(r)}{r}$, and
\begin{align}
u(p)= \frac{1}{n\omega_nr^{n- 1}}\int_{b= r}u|\nabla b|, \quad \quad \quad \quad \forall r> 0. \nonumber 
\end{align}
}
\end{lemma}

\begin{remark}
{The formula $I'_{u}(r)= 2\frac{D_{u}(r)}{r}$ firstly appeared in \cite{CM-JDG-97}, for reader's convenience we include its proof here.
}
\end{remark}

\pf
{Firstly, we note $\frac{\partial}{\partial r}= \frac{\nabla b}{|\nabla b|^2}$, assume the volume element of $b^{-1}(r)$ is $J(b^{-1}(r))$, then we have 
\begin{align}
I'(r)&= (1- n)r^{-n}\int_{b= r}u^2|\nabla b|+ r^{1- n}\int_{b= r}2u\langle \nabla u, \frac{\nabla b}{|\nabla b|}\rangle \nonumber \\
&\quad \quad + r^{1- n}\int_{b= r} u^2\Big(\frac{\nabla_{\nabla b}|\nabla b|}{|\nabla b|^2}+ \frac{\nabla_{\nabla b}J(b^{-1}(r))}{|\nabla b| J(b^{-1}(r))}\Big)dx \nonumber \\
&= 2\frac{D_u(r)}{r}+ (1- n)r^{-n}\int_{b= r}u^2|\nabla b|+ r^{1- n}\int_{b= r} \frac{u^2}{|\nabla b|}\Delta b \nonumber 
\end{align}
where we use the Divergence Theorem, $u$ is harmonic and the fact $\Big(\frac{\nabla_{\nabla b}|\nabla b|}{|\nabla b|}+ \frac{\nabla_{\nabla b}J(b^{-1}(r))}{J(b^{-1}(r))}\Big)= \Delta b$. From $\Delta (b^{2- n})= 0$, we have $\Delta b= \frac{n- 1}{b}|\nabla b|^2$, plug into the above equation, we get $I'_{u}(r)= 2\frac{D_{u}(r)}{r}$.

Similarly, we have 
\begin{align}
&\quad \quad \frac{d}{dr}\Big(r^{1- n}\int_{b= r}u|\nabla b|\Big) \nonumber \\
&= (1- n)r^{-n}\int_{b= r}u|\nabla b|+ r^{1- n}\int_{b= r} \langle \nabla u, \frac{\nabla b}{|\nabla b|}\rangle \nonumber \\
&\quad \quad + r^{1- n}\int_{b= r} u\Big(\frac{\nabla_{\nabla b}|\nabla b|}{|\nabla b|^2}+ \frac{\nabla_{\nabla b}J(b^{-1}(r))}{|\nabla b| J(b^{-1}(r))}\Big)dx \nonumber \\
&= r^{1- n}\int_{b\leq r}\Delta u+ (1- n)r^{-n}\int_{b= r}u|\nabla b|+ r^{1- n}\int_{b= r} \frac{u}{|\nabla b|}\Delta b \nonumber \\
&= (1- n)r^{-n}\int_{b= r}u|\nabla b|+ r^{1- n}\int_{b= r} \frac{u}{|\nabla b|}\Delta b= 0. \nonumber 
\end{align}
Finally, we get
\begin{align}
r^{1- n}\int_{b= r}u|\nabla b|= \lim_{r\rightarrow 0}r^{1- n}\int_{b= r}u|\nabla b|= u(p)\lim_{r\rightarrow 0}r^{1- n}\int_{b= r}|\nabla b|= u(p)n\omega_n . \nonumber 
\end{align}
}
\qed

\begin{lemma}\label{lem freq at the origin}
{If $u(p)= 0$, then $\mathscr{F}_u(0)= \frac{1}{2}\min\{m\geq 0: \frac{d^m}{dr^m} I_u(r)|_{r= 0}\neq 0\}\in \mathbb{Z}^+$.
}
\end{lemma}

\pf
{From Lemma \ref{lem derivative of I(r)} and by induction, it is straightforward to get 
\begin{align}
\frac{d^{m}}{dr^m}I(r)= r^{1- n}\int_{b= r} |\nabla b|\frac{\partial^m}{\partial r^m}(u^2), \quad \quad \quad \quad \forall m\geq 0. \nonumber 
\end{align}
Let $u\equiv 1$ in Lemma \ref{lem derivative of I(r)}, we get $r^{1- n}\int_{b= r}|\nabla b|= n\omega_n$. Hence
\begin{align}
\lim_{r\rightarrow 0}\frac{d^{m}}{dr^m}I(r)= n\omega_n\cdot \frac{\partial^m}{\partial r^m}(u^2)\big|_{r= 0} , \quad \quad \quad \quad \forall m\geq 0. \nonumber 
\end{align}

Let $j= \max\{s: \frac{\partial^t}{\partial r^t}u\big|_{r= 0}= 0, \ 0\leq t\leq s\}$ , then $\frac{\partial^t}{\partial r^t}(u^2)\big|_{r= 0}= 0$ for all $0\leq t\leq (2j+ 1)$, and $\frac{\partial^{2j+ 2}}{\partial r^{2j+ 2}}(u^2)\big|_{r= 0}\neq 0$. From Lemma \ref{lem derivative of I(r)}, we know that $D(r)= \frac{1}{2}r\cdot I'(r)$, hence $D^{(2j+ 2)}(r)\big|_{r= 0}= \frac{2j+ 2}{2}I^{(2j+ 2)}(r)\big|_{r= 0}$. By the L'Hospital rule, we have $\mathscr{F}(0)= \frac{D^{(2j+ 2)}(0)}{I^{(2j+ 2)}(0)}= j+ 1$.
}
\qed

\begin{cor}\label{cor local freq's increasing property}
{For any non-constant harmonic function $u$ defined on $B_{r}(p)\subseteq M^n$, there exists $C(p, r, M^n, u)< \infty$ such that $\displaystyle \sup\limits_{s\in [0, r]}\mathscr{F}_u(s)\leq C(p, r, M^n, u)$.
}
\end{cor}

\pf
{Note if $\mathscr{F}_u(s)= \infty$ for any $s> 0$, then $u\big|_{b(x)\leq s}= 0$. From the unique continuation property of harmonic functions, we know that $u\equiv 0$, it is the contradiction. Now the conclusion follows from Lemma \ref{lem freq at the origin} and the continuity of $\mathscr{F}_u(s)$ with respect to $s$.
}
\qed

The following lemma was proved in \cite{Li-Large}.
\begin{lemma}\label{lem Li's limit of subharmonic}
{Let $M^n$ be a complete manifold with $Rc\geq 0$. Suppose $f$ is a bounded subharmonic function defined on $M^n$, then $\lim\limits_{r\rightarrow \infty}\fint_{B_r(p)} f= \sup\limits_{M^n} f$.
}
\end{lemma}\qed

\begin{lemma}\label{lem est of I' and D}
{For any $\epsilon_0\in (0, 1)$, there exists $r_0> 0$ such that for any $r\geq (1-\epsilon_0)^{\frac{1}{2- n}}\mathrm{V}_M^{\frac{1}{n- 2}} r_0$, and harmonic function $u$ on $M^n$,
\begin{align}
(1- \epsilon_0)^{\frac{n}{2- n}}\fint_{\rho\leq (1- \epsilon_0)^{\frac{1}{2- n}}\mathrm{V}_M^{\frac{1}{2- n}} r} |\nabla u|^2\leq \frac{D(r)}{\omega_n\mathrm{V}_M^{\frac{2}{2- n}}r^2}\leq (1+ \epsilon_0)^{\frac{n+ 1}{n- 2}}\fint_{\rho\leq (1+\epsilon_0)^{\frac{1}{n- 2}}\mathrm{V}_M^{\frac{1}{2- n}} r} |\nabla u|^2 .\nonumber 
\end{align}
}
\end{lemma}

\pf
{From Lemma \ref{lem set inclusion between b and rho}, for $\epsilon_0> 0$, there is $r_0> 0$ such that for any $r\geq (1-\epsilon_0)^{\frac{1}{2- n}}\mathrm{V}_M^{\frac{1}{n- 2}} r_0$, we have
\begin{align}
D(r)&= r^{2- n}\int_{b\leq r}|\nabla u|^2\geq r^{2- n}\int_{\rho\leq (1- \epsilon_0)^{\frac{1}{2- n}}\mathrm{V}_M^{\frac{1}{2- n}} r} |\nabla u|^2 \nonumber \\
&\geq  \frac{\mathrm{Vol}(\rho\leq (1- \epsilon_0)^{\frac{1}{2- n}}\mathrm{V}_M^{\frac{1}{2- n}}r)}{r^{n- 2}}\fint_{\rho\leq (1-\epsilon_0)^{\frac{1}{2- n}}\mathrm{V}_M^{\frac{1}{2- n}} r} |\nabla u|^2\nonumber \\
&\geq (1- \epsilon_0)^{\frac{n}{2- n}}\omega_n\mathrm{V}_M^{\frac{2}{2- n}}r^2\fint_{\rho\leq (1- \epsilon_0)^{\frac{1}{2- n}}\mathrm{V}_M^{\frac{1}{2- n}} r} |\nabla u|^2 , \nonumber 
\end{align}
and 
\begin{align}
D(r)&= r^{2- n}\int_{b\leq r}|\nabla u|^2\leq r^{2- n}\int_{\rho\leq (1+\epsilon_0)^{\frac{1}{n- 2}}\mathrm{V}_M^{\frac{1}{2- n}} r} |\nabla u|^2 \nonumber \\
&\leq  \frac{\mathrm{Vol}(\rho\leq (1+ \epsilon_0)^{\frac{1}{n- 2}}\mathrm{V}_M^{\frac{1}{2- n}}r)}{r^{n- 2}}\fint_{\rho\leq (1+\epsilon_0)^{\frac{1}{n- 2}}\mathrm{V}_M^{\frac{1}{2- n}} r} |\nabla u|^2\nonumber \\
&\leq (1+ \epsilon_0)^{\frac{n+ 1}{n- 2}}\omega_n\mathrm{V}_M^{\frac{2}{2- n}}r^2\fint_{\rho\leq (1+\epsilon_0)^{\frac{1}{n- 2}}\mathrm{V}_M^{\frac{1}{2- n}} r} |\nabla u|^2 .\nonumber 
\end{align}
}
\qed

\begin{prop}\label{prop upper bound of freq for LGHF}
{For linear growth harmonic function $u$ on $M^n$, $\lim\limits_{s\rightarrow \infty}\mathscr{F}_u(s)= 1$.
}
\end{prop}

\pf
{Without loss of generality, we can assume $\sup\limits_{M^n}|\nabla u|= 1$. By Bochner formula and $u$ is harmonic, one gets that $|\nabla u|^2$ is a bounded subharmonic function. From Lemma \ref{lem est of I' and D}, for any $\epsilon_0\in (0, 1)$, there exists $r_0> 0$ such that if $r\geq r_0$, we have 
\begin{align}
\psi_1(\epsilon_0)\varphi\Big(\psi_2(\epsilon_0)\mathrm{V}_M^{\frac{1}{2- n}} r\Big)\leq \frac{D(r)}{\mathrm{V}_M^{\frac{2}{2- n}}r^2}\leq \psi_3(\epsilon_0)\varphi\Big(\psi_4(\epsilon_0)\mathrm{V}_M^{\frac{1}{2- n}} r\Big) , \nonumber  
\end{align}
where $\varphi(s)= \fint_{\rho\leq s} |\nabla u|^2$ and $\psi_i(\cdot)$ are functions satisfying $\lim\limits_{t\rightarrow 0}\psi_i(t)= 1, i= 1, 2, \cdots$.

The above inequality implies that for any $s\in [r_0, r]$,
\begin{align}
\psi_5(\epsilon_0)\frac{s^2}{r^2}\cdot \frac{\varphi\Big(\psi_2(\epsilon_0)\mathrm{V}_M^{\frac{1}{2- n}} s\Big)}{\varphi\Big(\psi_4(\epsilon_0)\mathrm{V}_M^{\frac{1}{2- n}} r\Big)}\leq \frac{D(s)}{D(r)}\leq \psi_6(\epsilon_0)\frac{s^2}{r^2}\cdot \frac{\varphi\Big(\psi_4(\epsilon_0)\mathrm{V}_M^{\frac{1}{2- n}} s\Big)}{\varphi\Big(\psi_2(\epsilon_0)\mathrm{V}_M^{\frac{1}{2- n}} r\Big)} .\nonumber 
\end{align}
Combining Lemma \ref{lem Li's limit of subharmonic}, if $r_0$ is big enough, we have
\begin{align}
\psi_5(\epsilon_0)\frac{s^2}{r^2}(1- \epsilon_0)\leq \frac{D(s)}{D(r)}\leq \psi_6(\epsilon_0)\frac{s^2}{r^2}\cdot (1+ \epsilon_0), \nonumber 
\end{align}
which implies $\displaystyle \psi_7(\epsilon_0)\frac{s^2}{r^2}\leq \frac{D(s)}{D(r)}\leq \psi_8(\epsilon_0)\frac{s^2}{r^2}$. From Lemma \ref{lem derivative of I(r)}, 
\begin{align}
\mathscr{F}(r)= \frac{D(r)}{I(r_0)+ \int_{r_0}^r I'(s)ds}= \frac{D(r)}{I(r_0)+ \int_{r_0}^r\frac{2D(s)}{s}ds}= \frac{1}{\frac{I(r_0)}{D(r)}+ 2\int_{r_0}^r\frac{1}{s}\frac{D(s)}{D(r)}ds} .\nonumber 
\end{align}

Note $I(r_0)\geq 0$, if we assume $r\geq 2r_0$, we get
\begin{align}
\mathscr{F}(r)&\leq \frac{1}{2\int_{r_0}^r\frac{1}{s}\frac{D(s)}{D(r)}ds}\leq \frac{1}{2\psi_7(\epsilon_0)\int_{r_0}^r \frac{s}{r^2}ds}= \frac{\psi_9(\epsilon_0)}{1- \frac{r_0^2}{r^2}} ,\nonumber 
\end{align}
which implies $\lim\limits_{r\rightarrow \infty}\mathscr{F}(r)\leq \psi_9(\epsilon_0)$. Letting $\epsilon_0\rightarrow 0$, we get $\lim\limits_{r\rightarrow \infty}\mathscr{F}(r)\leq 1$.

On the other hand, similar to the above, we get
\begin{align}
\mathscr{F}(r)\geq \frac{1}{\frac{I(r_0)}{D(r)}+ 2\psi_8(\epsilon_0)\int_{r_0}^r\frac{s}{r^2}ds}\geq \frac{1}{\frac{I(r_0)}{\psi_1(\epsilon_0)\varphi\Big(\psi_2(\epsilon_0)\mathrm{V}_M^{\frac{1}{2- n}} r\Big)\mathrm{V}_M^{\frac{2}{2- n}}r^2}+ \psi_8(\epsilon_0)\big[1- \frac{r_0^2}{r^2}\big]} ,\nonumber 
\end{align}
which implies $\lim\limits_{r\rightarrow \infty}\mathscr{F}(r)\geq 1$. The conclusion follows.
}
\qed

Recall S.-Y. Cheng \cite{Cheng-sub} proved that any non-constant polynomial growth harmonic function is at least linear growth. From \cite[Theorem $1$]{CCM}, we know that on $M^n$ with $Rc\geq 0$, the existence of non-constant linear growth harmonic function implies, that any tangent cone at infinity of $M^n$, $M_{\infty}$, will split isometrically as $\mathbb{R}\times N$. The following rigidity result links the geometric structure of $M^n$ with the frequency upper bound of global harmonic functions.
\begin{cor}\label{cor constant freq harmonic implies Rn}
{For harmonic function $u$ with $u(p)= 0$, if $\sup\limits_{s\geq 0}\mathscr{F}_u(s)\leq 1$, then $M^n$ is isometric to $\mathbb{R}^n$. 
}
\end{cor}

\pf
{From Lemma \ref{lem poly growth by frequency} and $\mathscr{F}(s)\leq 1$, we know that $u$ is linear growth. Hence we can assume that $\sup\limits_{M^n}|\nabla u|= 1$ by rescaling. From Lemma \ref{lem derivative of I(r)}, we have $\frac{I'}{I}(r)= \frac{2}{r}\mathscr{F}(r)$, which implies $I(r)= I(1)e^{\int_1^r \frac{2\mathscr{F}(s)}{s}ds}$ and
\begin{align}
D(r)= I(r)\mathscr{F}(r)= I(1)\mathscr{F}(r)e^{\int_1^r \frac{2\mathscr{F}(s)}{s}ds}. \nonumber 
\end{align}

By Lemma \ref{lem set inclusion between b and rho} and Proposition \ref{prop upper bound of freq for LGHF}, for any $\epsilon_0> 0$ the followings hold:
\begin{align}
I(1)&= \lim_{r\rightarrow\infty}\frac{D(r)}{\mathscr{F}(r)e^{\int_1^r \frac{2\mathscr{F}(s)}{s}ds}}= \lim_{r\rightarrow\infty}\frac{r^{2- n}\int_{b\leq r}|\nabla u|^2}{\mathscr{F}(r)e^{\int_1^r \frac{2\mathscr{F}(s)}{s}ds}} \nonumber \\
&\geq \lim_{r\rightarrow\infty} \frac{r^2}{\mathscr{F}(r)e^{\int_1^r \frac{2\mathscr{F}(s)}{s}ds}}\cdot \frac{\int_{\rho\leq (1- \epsilon_0)^{\frac{1}{2- n}}r\mathrm{V}_M^{\frac{1}{2- n}}}|\nabla u|^2}{r^n} \nonumber \\
&= \varliminf_{r\rightarrow\infty} \frac{r^2}{\mathscr{F}(r)e^{\int_1^r \frac{2\mathscr{F}(s)}{s}ds}}\cdot \lim_{r\rightarrow\infty} \frac{V\Big(\rho\leq r\mathrm{V}_M^{\frac{1}{2-n}}(1- \epsilon_0)^{\frac{1}{2- n}}\Big)}{r^n}\sup_{M^n}|\nabla u|^2 \nonumber \\
&= \omega_n\mathrm{V}_M^{\frac{2}{2- n}}(1- \epsilon_0)^{\frac{n}{2- n}} . \nonumber 
\end{align}

On the other hand, by $u(p)= 0$ and Lemma \ref{lem freq at the origin}, we know $\lim\limits_{r\rightarrow 0}\mathscr{F}(r)\geq 1$, hence $\mathscr{F}(0)= 1$ and
\begin{align}
I(1)&= \lim_{r\rightarrow 0}\frac{D(r)}{\mathscr{F}(r)e^{\int_1^r \frac{2\mathscr{F}(s)}{s}ds}}\leq \lim_{r\rightarrow 0}\frac{D(r)}{r^2}\cdot \frac{r^2}{e^{\int_1^r \frac{2\mathscr{F}(s)}{s}ds}}\nonumber \\
&\leq \lim_{r\rightarrow 0}\frac{\int_{b\leq r} |\nabla u|^2}{r^n}\cdot \varlimsup_{r\rightarrow 0}\frac{r^2}{e^{\int_1^r \frac{2\mathscr{F}(s)}{s}ds}}= \omega_n|\nabla u|^2(p). \nonumber 
\end{align}

Combining the above, we get 
\begin{align}
|\nabla u|^2(p)\geq (1- \epsilon_0)^{\frac{n}{2- n}}\mathrm{V}_M^{\frac{2}{2- n}}\geq (1- \epsilon_0)^{\frac{n}{2- n}}. \nonumber 
\end{align}
Letting $\epsilon_0\rightarrow 0$ in the above, we have $|\nabla u|^2(p)\geq  1$. From the strong maximum principle and $\Delta |\nabla u|^2\geq 0$, we get $|\nabla u|^2\equiv 1$. This implies $\mathrm{V}_M= 1$, hence $M^n$ is isometric to $\mathbb{R}^n$.
}
\qed

Now we present the proof of Theorem \ref{thm uniform bound of freq}.

\pf[of Theorem \ref{thm uniform bound of freq}]
{From \cite{LT-linear}, $\mathrm{dim}\mathcal{H}_1(M^n)< \infty$, where $\mathcal{H}_1(M^n)$ is the space of linear growth harmonic functions on $M^n$. We can assume $\mathcal{H}'(M)= \{f: f(p)= 0, f\in \mathcal{H}_1(M)\}$, and define the inner product as
\begin{align}
\langle f, g\rangle= \lim_{r\rightarrow \infty}\fint_{B_r(p)}\nabla f\cdot \nabla g , \nonumber 
\end{align}
which is well-defined by Lemma \ref{lem Li's limit of subharmonic}.

Assume $\{u_i\}_{i= 1}^m$ be an orthonormal basis of $\mathcal{H}'(M)$, to prove the conclusion, we only need to show that $\mathscr{F}_u(r)\leq C(M^n)$, where $\displaystyle u= \sum_{i= 1}^m a_iu_i$ with $\displaystyle \sum_{i= 1}^m |a_i|^2= 1$.

For any $\epsilon_0> 0$, there is $R_0> 0$ such that if $t\geq R_0$, we have
\begin{align}
&\Big|\fint_{B_{\psi_k(\epsilon_0)\mathrm{V}_M^{\frac{1}{2- n}}t}(p)}\nabla u_i\cdot \nabla u_j- \delta_{ij}\Big|\leq \epsilon_0, \quad \quad \quad \quad 1\leq i, j\leq m, 1\leq k\leq 4 ,\nonumber  \\
&\Big|\fint_{B_{\psi_k(\epsilon_0)\mathrm{V}_M^{\frac{1}{2- n}}t}(p)}|\nabla u|^2- 1\Big|\leq \sum_{i, j= 1}^m |a_ia_j|\cdot \Big|\fint_{B_{\psi_k(\epsilon_0)\mathrm{V}_M^{\frac{1}{2- n}}t}(p)} \nabla u_i\cdot \nabla u_j- \delta_{ij}\Big| \nonumber \\
&\leq \epsilon_0\Big(\sum_{i= 1}^m |a_i|\Big)^2 \leq m\epsilon_0. \nonumber 
\end{align}

Then from Lemma \ref{lem est of I' and D}, for $r\geq s\geq R_0$, we get
\begin{align}
\frac{D_u(s)}{D_u(r)}\geq \frac{\psi_1(\epsilon_0)s^2}{\psi_3(\epsilon_0)r^2}\frac{\fint_{B_{\psi_2(\epsilon_0)\mathrm{V}_M^{\frac{1}{2- n}}s}(p)}|\nabla u|^2}{\fint_{B_{\psi_4(\epsilon_0)\mathrm{V}_M^{\frac{1}{2- n}}r}(p)}|\nabla u|^2}\geq \frac{1- m\epsilon_0}{1+ m\epsilon_0}\psi(\epsilon_0)\frac{s^2}{r^2} . \nonumber 
\end{align}
And if $r\geq 2R_0$ we get 
\begin{align}
\mathscr{F}(r)&= \frac{1}{\frac{I(R_0)}{D(r)}+ 2\int_{R_0}^r\frac{1}{s}\frac{D(s)}{D(r)}ds}\leq \frac{1}{2\int_{R_0}^r\frac{1}{s}\frac{D(s)}{D(r)}ds}\leq \frac{1}{2\psi(\epsilon_0)\int_{R_0}^r \frac{s}{r^2}ds}= \frac{\psi(\epsilon_0)}{1- \frac{(R_0)^2}{r^2}} \nonumber \\
&\leq \frac{4}{3}\psi(\epsilon_0).\nonumber
\end{align}

For $r\leq 2R_0$, we claim that there is some $C> 0$ such that $\mathscr{F}_u(r)\leq C$ for any $u\in \mathcal{H}'(M^n)$. By contradiction, if there is $v_i\in \mathcal{H}'(M^n)$ and $s_i\leq 2R_0$ such that $\lim\limits_{i\rightarrow \infty}\mathscr{F}_{v_i}(s_i)= \infty$. Without loss of generality, we can rescale $v_i$ such that $\max\limits_{x\in B_{2R_0}(p)}|\nabla v_i(x)|= 1$. From the compactness theorem for harmonic functions, we get $\lim\limits_{i\rightarrow \infty}s_i= s_{\infty}\leq 2R_0$ and $v_i\rightarrow v_\infty\in \mathcal{H}_1(M^n)$, also
\begin{align}
\mathscr{F}_{v_\infty}(s_\infty)= \lim_{i\rightarrow\infty}\mathscr{F}_{v_i}(s_i)= \infty, \quad \quad \quad v_\infty(p)= 0, \quad \quad \quad \max_{x\in B_{2R_0}(p)}|\nabla v_\infty(x)|= 1 ,\nonumber 
\end{align}
which is the contradiction to Corollary \ref{cor local freq's increasing property}. The conclusion follows from the above argument.
}
\qed

\section{Polynomial growth harmonic functions on Perelman's manifold}

From \cite{Ding} (also see \cite{Xu}), we know the existence of non-constant PGHF on maximal volume growth manifolds with unique tangent cone at infinity. In this section, applying the results of \cite{Xu}, we show the existence of non-constant PGHF on some maximal volume growth manifolds with different tangent cones at infinity, which were constructed firstly by Perelman \cite{Perelman-example}.

A metric space $(M_{\infty}, p_{\infty}, \rho_{\infty})$ is a \textbf{tangent cone at infinity} of $M^n$ if it is a Gromov-Hausdorff limit of a sequence of rescaled manifolds $(M^n, p, r_j^{-2}g)$, where $r_j\rightarrow \infty$. By Gromov's compactness theorem, \cite{Gromov}, any sequence $r_j\rightarrow \infty$, has a subsequence, also denoted as $r_j\rightarrow \infty$, such that the rescaled manifolds $(M^n, p, r_j^{-2}g)$ converge to some tangent cone at infinity $M_{\infty}$ in the Gromov-Hausdorff sense. 
From Cheeger-Colding's theory of Ricci limit spaces, there exists a self-adjoint Laplace operator $\Delta_{(C(X), \nu)}$ on $(C(X), \nu)\in \mathscr{M}(M)$, where $\nu$ is a measure. And $\nu$ induces a natural measure $\nu_{-1}$ on $X$, which yields the existence of a self-adjoint positive Laplace operator $\Delta_{(X, \nu_{-1})}$ on $(X, \nu_{-1})$.

On a metric cone $(C(X), dr^2+ r^2dX)$, the measure $\nu$ is called \textbf{conic measure of power $\kappa$},  and $\kappa$ is a positive constant denoted as $\mathtt{p}(\nu)$, if for any $\Omega\subset\subset C(X)$, 
\begin{align}
\nu (\Omega)= \int_0^{\infty} r^{\kappa- 1}dr\int_{X} \chi(\Omega_r) d\nu_{-1} \nonumber 
\end{align}
where $\Omega_r= \{z| z\in \Omega, r(z)= r\}$, $\chi(\cdot)$ is the characteristic function on $C(X)$.

And we also define $\mathscr{S}(M)$ the spectrum at infinity of $(M^n, g)$ and $\mathscr{D}(M)$ the degree spectrum at infinity of $(M^n, g)$:
\begin{align}
\mathscr{S}(M)&\vcentcolon= \{\lambda|\ \lambda= \lambda_j(X, \nu_{-1})\ for \ some\ positive\ interger\  j\ and \ (C(X), \nu)\in \mathscr{M}(M)\} \nonumber \\
\mathscr{D}(M)&\vcentcolon= \{\alpha\geq 0|\ \alpha\big(\kappa+ \alpha- 2\big)= \lambda\ for \ some\ \lambda= \lambda_j(X, \nu_{-1})\in \mathscr{S}(M)\  and \ \kappa= \mathtt{p}(\nu)\} \nonumber
\end{align}

Recall the following theorem \cite[Theorem $4.3$]{Xu}.
\begin{theorem}\label{thm general existence of poly growth harmonic function}
{Let $(M^n, g)$ be a complete manifold with nonnegative Ricci curvature, assume that every tangent cone at infinity of $M^n$ with renormalized limit measure is a metric cone $C(X)$ with conic measure of power $\kappa\geq 2$, and $\mathcal{H}^{1}(X)> 0$. If there exists $d\notin \mathscr{D}(M)$ and $d> \inf \{\alpha|\ \alpha\in \mathscr{D}(M), \alpha\neq 0\}$, then $dim\big(\mathscr{H}_d(M)\big)\geq 2$.
}
\end{theorem}\qed

The following lemma was implied by Perelman's work \cite{Perelman} (also see \cite{AMW}).
\begin{lemma}\label{lem sharp angle cone can be smoothed off}
{For any Riemannian manifold $(M^4, g)$ with $Rc(g)> 0$, and assume the boundary $\partial M^4$ is isometric to round sphere $\mathbb{S}^{3}$, whose sectional curvature is equal to $r_1^2> 0$, furthermore its principal curvatures $\lambda(\partial M^4)$ satisfy $\displaystyle \sup \big[\lambda(\partial M^4)\big]^2\leq c_1^2 r_1^2$, where $c_1\in (0, 1)$ is some constant. Then the differential manifold $M^4\bigcup_{\mathbb{S}^{3}} (\mathbb{S}^4- \mathcal{B}^4)$ carries a Riemannian metric $\hat{g}$ with $Rc(\hat{g})> 0$, where $\mathcal{B}^4$ is diffeomorphic to a $4$-dim Euclidean ball.
}
\end{lemma}\qed

To show the existence of non-constant polynomial growth harmonic functions on Perelman's example manifolds, we recall the construction of  Perelman's example  \cite{Perelman-example} for completeness. For these manifolds with $Rc\geq 0$ has maximal volume growth, the tangent cones at infinity are not unique.

Consider the metric $ds^2= dt^2+ \hat{g}_{\mathbb{S}^3}$ defined on $\mathbb{R}^4$, where $t\geq 0$, $\hat{g}_{\mathbb{S}^3}= A^2(t)dX^2+ B^2(t)dY^2+ C^2(t)dZ^2$, $X, Y, Z$ are vector fields on $\mathbb{S}^3$, $T= \frac{\partial}{\partial t}$, and $X, Y, Z, T$ satisfy  
\begin{align}
[X, T]= [Y, T]= [Z, T]= 0, \quad \quad [X, Y]= 2Z, \quad \quad [Y, Z]= 2X, \quad \quad [Z, X]= 2Y \nonumber 
\end{align}

Now we take 
\begin{align}
A(t)&= \frac{t}{10}\big[1+ \phi(t)\sin(\ln\ln t)\big] , \quad \quad \quad B(t)= \frac{t}{10}\cdot \frac{1}{1+ \phi(t)\sin(\ln\ln t)} ,\nonumber \\
C(t)&= \frac{t}{10}\big[1- \gamma(t)\big] ,\nonumber
\end{align}
where $\phi(t)$ is a smooth function such that
\begin{align}
&\phi(t)= 0 \ , \quad \quad \forall t\in [0, t_0], \nonumber \\
&\phi(t)> 0 \ , \quad \quad \forall t> t_0 ,\nonumber \\
&0\leq \phi'(t)\leq \frac{1}{t^2} \quad \quad and \quad \quad |\phi''(t)|\leq \frac{1}{t^3}, \nonumber
\end{align}
and $\gamma(t)$ is a smooth function defined on $[\frac{t_0}{2}, \infty)$ such that 
\begin{align}
&\gamma(\frac{t_0}{2})= 0, \ \gamma'(t)> 0\ , \ \gamma''(t)> 0 \ ,\quad\quad \forall t\in (\frac{t_0}{2}, t_0) ,\nonumber \\
&\gamma'(t)= \frac{1}{t\ln^{\frac{3}{2}}t} \ , \quad \quad \forall t> t_0 , \nonumber 
\end{align}
We consider the manifold $(N^4, ds^2)= \big(\mathbb{R}^4- B_{\frac{t_0}{2}}(0), ds^2\big)$, where $B_{\frac{t_0}{2}}(0)$ is the metric ball with respect to the metric $ds^2$. One can check that $Rc(N^4, ds^2)> 0$ from the above properties of $g$.

Note near $\partial N^4$, $A= B= C= \frac{t}{10}$, the principal curvature $\lambda(\partial N^4)$ is $-\frac{1}{t}$. Because $\partial N^4$ is round sphere, the intrinsic sectional curvatures $\widetilde{Rm}(\partial N^4)$ of $\partial N^4$ are $\frac{10^2}{t^2}$.

Hence $\big[\lambda(\partial N^4)\big]^2= \frac{1}{100}\widetilde{Rm}(\partial N^4)$. Using Lemma \ref{lem sharp angle cone can be smoothed off}, one gets the differential manifold $M^4= N^4\bigcup_{\mathbb{S}^3} (\mathbb{S}^4- \mathcal{B}^4)$, which admits a Riemannian metric $g$ with $Rc(g)> 0$, and $g= ds^2$ outside of a compact subset of $M^4$. 

We can define $\displaystyle \phi_\infty= \lim_{t\rightarrow \infty}\phi(t)> 0$ and $\displaystyle\gamma_\infty= \lim_{t\rightarrow \infty}\gamma(t)> 0$. Then $\displaystyle \mathcal{M}_\infty$ is a family of $\big(C(X), dt^2+ t^2dX\big)$, where $X= (\mathbb{S}^3, g_X)\Big\}$ with
\begin{align}
g_X= \Big(\frac{1+ \phi_\infty\delta}{10}\Big)^2dX^2+ \Big(\frac{1}{10}\cdot \frac{1}{1+ \phi_\infty\delta}\Big)^2dY^2+ \Big(\frac{1- \gamma_\infty}{10}\Big)^2dZ^2, \quad \quad \forall \delta\in [-1, 1].  \nonumber 
\end{align}
And the conic measure is of power $4$ for each $\delta$.

If we choose $\phi_\infty, \gamma_\infty$ small enough, then $\mathscr{S}(M^4)$ is a set close to the spectrum of $\mathbb{S}^3$, which implies the existence of $d\notin \mathscr{D}(M^4)$ satisfying $\displaystyle d> \inf \mathscr{D}(M^4)$. From \cite[Theorem $4.3$]{Xu}, one gets the existence of polynomial growth harmonic functions on Perelman's example manifolds.

\section*{Acknowledgments}
The author thank the anonymous referee for detailed suggestion on the former version of the paper, which improves the the presentation of the paper. We thank Christina Sormani for her interest and comments.



\end{document}